\documentclass{article}
\usepackage{graphicx}
\usepackage{amssymb,amsfonts,amsthm}
\usepackage{amsmath,amsthm,amssymb}
\usepackage{amsfonts}

\newtheorem{theorem}{Theorem}[section]
\newtheorem{lemma}[theorem]{Lemma}

\theoremstyle{definition}

\newtheorem{rk}[theorem]{Remark}

\newcounter{ppp}

\newcommand{\la}{\langle}
\newcommand{\ra}{\rangle}

\begin{document}

\renewcommand{\theequation}{\thesection.\arabic{equation}}

\title{On products of T-ideals in free algebras and free group algebras}
 \author{ A. Yu. Olshanskii \thanks{\large The
author was supported in part by the NSF grants DMS
0455881, DMS 0700811, and by the Russian Fund for Basic Research grant
05-01-00895.}}
\date{}
\maketitle

\begin{abstract} Let $F$ be a field and $A$ a free
associative $F$-algebra or a group algebra of a free group with an infinite set $X$ of generators. We find a necessary and sufficient condition
for the inclusion $I'\subseteq I$, where $I=I_1\dots I_k$  and $I'=I'_1\dots I'_l$ are any products of $T$-ideals in $A$. 
A canonical reformulation in terms of products of group representation varieties
answers a question posed in 1986 \cite{O}.   

\end{abstract}

\section{Introduction}

Let $F\la X \ra =F\la x_1,x_2,\dots\ra$ be a free associative algebra over a field $F$ with $1$
and with a countable set of free non-commuting  generators $X=\{x_1,x_2,\dots\}$.
Recall that an ideal $I$ of $F\la X \ra$ is called a $T$-{\it ideal} if $\mu(I)\subseteq I$ for
every endomorphism $\mu$ of $F\la X \ra$. As usual, a product $IJ$ of ideals $I$ and $J$
is the set of finite sums $\sum a_i b_i$ where $a_i\in I, b_i\in J$. We prove

\begin{theorem} \label{theor} Let $V$ and $V'$ be products of $T$-ideals in $F\la X\ra$: 
$V=V_1\dots V_k$, and $V'=V'_1\dots V'_l$. Then we have $V'\subseteq V$ if and only if there exists a
$T$-ideal $I$ of $F\la X\ra$ with two $T$-ideal factorizations 

$$I=I_1\dots I_k=I'_1\dots I'_l$$
such that $$I_1\subseteq V_1,\dots, I_k\subseteq V_k \;\; and \;\; V'_1\subseteq I'_1,\dots,V'_l\subseteq I'_l.$$
\end{theorem}

The elements of $F\la X\ra$ are called (non-commutative) {\it polynomials} in this paper.
 Recall that an $F$-algebra $R$ satisfies an identity $f(x_1,\dots,x_n)=0$ for a polynomial $f(x_1,\dots,x_n)=f$ if
 $f$ vanishes under every homomorphism $F\la X\ra\to R$. A {\it variety} of associative algebras
 is a class of algebras consisting of all associative algebras with $1$ that satisfy some set of identities. 
 For every algebra $R$, the set of the
 left-hand sides of all its identities $f=0$ is a $T$-ideal in $F\la X\ra$, and   it is well known that the set of 
 varieties is in Galois correspondence with the set of $T$-ideals of  $F\la X\ra$. 

Let  $\cal U$ and $\cal V$ be varieties of associative
algebras. We denote by $U$ and $V$, respectively, the corresponding 
$T$-ideals in $F\la X \ra$.
The {\it product} ${\cal U}{\cal V}$ of the varieties $\cal U$ and $\cal V$ is
defined by the $T$-ideal $UV$, i.e., an algebra $R$ belongs to ${\cal U}{\cal V}$ iff
it satisfies all the identities $uv=0$, where $u\in U$ and $v\in V$.

Thus, Theorem \ref{theor} admits an obvious reformulation in term of products
of varieties. Such a reformulation is a strict analog of the theorem on products
of group varieties and of the theorem on products of the varieties of Lie algebras
over an infinite field from \cite{O} and \cite{OS}, respectively.

In turn, the 
papers \cite{O} and \cite{OS} extend  the equality case   
${\cal V}_1\dots {\cal V}_k={\cal V'}_1\dots {\cal V'}_l$
considered earlier for group varieties in \cite{3N} and \cite{S}, and for varieties
of Lie algebras over an infinite field, in \cite{P} and \cite{B}. Indeed, "$=$" is
a special case of "$\subseteq$", and the
Neumanns--Shmel'kin theorem (the Parfenov--Bahturin theorem)  says that the
non-trivial group varieties (varieties of algebras Lie over an infinite field)
form a free monoid under the product operation. (Recall that the product
${\cal U}{\cal V}$ of two group varieties contains all groups $G$ having a normal
subgroup $N$ such that $N\in \cal U$ and $G/N\in \cal V$. The definition of
the product for the varieties of Lie algebras is similar, although such a multiplication
is not associative over finite fields \cite{B1}.)

It might seem that the associative algebra case differs from the group and Lie algebra
ones because the ideals of $F\la X\ra$ are usually not free. But instead of this,
all right (and left) ideals are free modules over $F\la X\ra$, and making use of this
fact, George Bergman and Jacques Lewin \cite{BL} proved that the monoid of non-zero ideals
of $F\la X\ra$ is free, and furthermore, the monoid of non-zero $T$-ideals of $F\la X\ra$
is also free.

The lemmas of the present papers are based on the Schreier--Lewin technique 
for firs (i.e., free ideal rings, see \cite{C}) and on the technique of triangular products of modules
invented by B.I.Plotkin (see \cite{PG}, \cite{PV}, \cite{K}, \cite{V}). To complete
the proof of Theorem \ref{theor}, we follow the outline from \cite{O}.

The triangular products were originally introduced and applied to the products
of varieties of group representations (see \cite{PV}, \cite{V}). Let $FG$ be 
the group algebra of a free group $G$ with
an infinite set of free generators $X=\{x_1,x_2,\dots\}$ over a field $F$, $f\in FG$ a 'polynomial',
and $H$ a group. Recall that an $FH$-module $M$ is said to satisfy the identity $f=0$ if
for every $a\in M$ and every homomorphism $\mu: G\to H$,
 we have $a(\bar\mu(f))=0$,
where $\bar\mu$ is the algebra homomorphism $FG\to FH$ induced by $\mu$.
A {\it  variety of group representations} is the class of group representations (or modules
over group algebras over $F$) satisfying a set of representation identities. 
A product $\cal UV$ of such varieties
is a variety consisting of all $FG$-modules $M$ having a submodule $N$ such that
$N\in\cal U$ and $M/N\in \cal V$. This multiplication is obviously associative,
and B.I.Plotkin \cite{PG} proved that the monoid of non-zero varieties of group
representations over any field $F$ is free. The formulation and the proof of the following theorem 
are similar to those for \ref{theor}. 

\begin{theorem}\label{theor'} Let $\cal V$ and $\cal V'$ be products of varieties of group
representations over a field: 
${\cal V}={\cal V}_1\dots {\cal V}_k$, and ${\cal V'}={\cal V'}_1\dots {\cal V'}_l$. 
Then we have $\cal V\subseteq \cal V'$ if and only if there exists a
variety $\cal W$ with two factorizations ${\cal W}={\cal W}_1\dots {\cal W}_k={\cal W'}_1\dots {\cal W'}_l$
such that ${\cal V}_1\subseteq {\cal W}_1,\dots, {\cal V}_k\subseteq {\cal W}_k$ and 
${\cal W'}_1\subseteq {\cal V'}_1,\dots,{\cal W'}_l\subseteq {\cal V'}_l$.
\end{theorem}

Theorem \ref{theor'} solves a problem raised in \cite{O}.  

\section{Free module bases of $T$-ideals}

We denote by $X^*$ the free monoid generated by $X=\{x_1,x_2,\dots\}$. Recall that a subset $S\subseteq X^*$ 
is called a {\it Schreier set} if with any monomial (word) it contains all its left factors
(prefixes).   

Let $I$ be a right ideal of $F\la X\ra$. Following Jacques Lewin \cite{L}, we call
a set $S$ of monomials  a {\it Schreier basis} for $F\la X\ra$
modulo $I$ if $S$ is a Schreier set 
  and the image of $S$ in $F\la X\ra/I$ is an $F$-basis of the quotient $F\la X\ra/I$.

We denote by $\phi$ and $\psi$ the endomorphisms of the algebra $F\la X\ra$ such that $\phi(x_i)=x_{i+1}$ for $i=1,2,\dots$
and $\psi (x_{i+1})=x_i$ for $i=1,2,\dots$, $\psi (x_1)=0$. Note that $\phi(I)\subseteq I$ and $\psi(I)\subseteq I$
if $I$ is a $T$-ideal of $F\la X\ra$.

 \begin{lemma}\label{schr}
For any $T$-ideal $I$ of $F\la X \ra$, there is a Schreier basis 
$S_I$ of $F\la X \ra$ modulo $I$ such that $\phi(S_I)\subseteq S_I$.
\end{lemma}

\proof 

Assume that the monoid $X^*$ is well ordered by degrees and lexicographically for the monomials of the same degree. 
 It is convenient
 to set $0<u$ for every $u\in X^*$. Then for $u_1, u_2, u\in X^*$, we have $\phi(u_1)\le \phi (u_2)$ iff $u_1\le u_2$,
 $\psi(u)<u$,
 and, provided $\psi (u_2)\ne 0$, we have $\psi(u_1)\le \psi(u_2)$ if $u_1\le u_2$.

Call an element of $X^*$ $I$-{\it reducible} if it lies in the sum of
$I$ and the span of the elements preceding it in the given ordering;
call it $I$-{\it irreducible} otherwise. Then the images in $F\la X\ra/I$ of
$I$-irreducible elements of $X^*$ form a basis. Denote the set of $I$-irreducible
elements of $X^*$ by $S_I$.

That the set $S_I$ is a Schreier set is
equivalent to saying that if a left factor of a word is $I$-reducible,
so is the whole word; but this is immediate (and only uses the definition of the
well-ordering and the fact
that  $I$  is a right ideal).

To see that $\phi(S_I)\subseteq S_I$,  consider any  $u\in S_I$  such that  $\phi(u)$  is reducible.  Thus,
$\phi(u)-r\in I$  where  $r$  is an $F$-linear combination of monomials  $< \phi(u)$.
Applying  $\psi$  to this relation, we see that  $I$  contains  $u-\psi(r)$,
so  $u$  is also reducible; so as the contrapositive, we have that if  $u$
is irreducible, so is  $\phi(u)$.

 \endproof

\begin{lemma}\label{basis} Every $T$-ideal $I$ of $F\la X\ra$ is a free (right) $F\la X\ra$-module
having an $F\la X\ra$-basis $B$ such that $\phi(B)\subseteq B$.
\end{lemma}
    
\proof It is well-known that $F\la X\ra$ is a free right (and left) ideal ring or {\it right (and left) fir},
i.e., every right ideal $I$ of $F\la X\ra$ is a free $F\la X\ra$-module. Furthermore, all non-zero 
differences of the form $e_{\alpha}=s x_i -\sum \lambda_j s_j$ form an $F\la X\ra$-basis $B$ of the right module $I$ (\cite{C}, Theorem VI.6.8).
Here $s$ is an arbitrary monomial from the Schreier $F$-basis $S_I$ of $F\la X\ra$ modulo $I$, $x_i$
is any element of $X$, and $\sum \lambda_j s_j$ is the unique linear combination of  vectors $s_j$
from $S_I$ such that $s x_i=\sum \lambda_j s_j$ ($mod\; I$).

If $e_{\alpha}\ne 0$, then $\phi(e_{\alpha})\in I \backslash \{0\}$ since $\phi$ is injective and
$I$ is a $T$-ideal. Then $\phi(e_{\alpha})=\phi(s)x_{i+1}-\sum \lambda_j\phi(s_j)$, where $\phi(s),\phi(s_j)\in S_I$
by Lemma \ref{schr}. Hence $\phi(e_{\alpha})$ has the same form and therefore belongs to the basis $B$.
The lemma is proved. \endproof

The following is a well known property of firs.

\begin{lemma}\label{incl} Let $U$, $I$ and $J$ be ideals of $F\la X\ra$ and $U\ne \{0\}$. Then $UI\subseteq UJ$ iff $I\subseteq J$.
\end{lemma} 

\proof The ideal $U$ has a non-empty basis $(e_{\alpha})$ as right $F\la X\ra$-module. Therefore a sum $\sum e_{\alpha}f_{\alpha}$
belongs to the product $UI$ (to $UJ$) iff all the polynomials $f_{\alpha}$ belong to $I$ (to $J$). This proves the statement of the lemma. \endproof  

\section{Annihilators and separators}

Let $M$ be a right $F\la X\ra$-module and $N$ a submodule. For a subset $Z\subseteq M$, its
annihilator modulo $N$ is

$$(N:Z)=\{a \mid a\in F\la X\ra, za\in N \;\; for\;\; all\;\; z\in Z\}$$
   The annihilator $(N:Z)$ is an ideal of $F\la X\ra$ if $Z$ is a submodule of $M$.

\begin{lemma}\label{ann} If $I$ and $J$ are $T$-ideals of $F\la X\ra$, then $(I:J)$ is also
a $T$-ideal in $F\la X\ra$.

\end{lemma} 
\proof We must show that $\mu(a)\in(I:J)$ for every $a\in(I:J)$ and every
endomorphism $\mu$ of $F\la X\ra$. Let $z=z(x_1,\dots,x_n)\in J$ and $a=a(x_1,\dots,x_n)$.
Since $J$ is a $T$-ideal, the polynomial $z(x_{n+1},\dots,x_{2n})$ also belongs to $J$,
and therefore $z(x_{n+1},\dots,x_{2n}) a(x_1,\dots,x_n)\in I$. Now we note that there is an endomorphism $\eta$
of $F\la X\ra$ such that $\eta(x_1)=\mu(x_1), \dots, \eta(x_n)=\mu(x_n),\; \eta(x_{n+1})=x_1,\dots,\eta(x_{2n})=x_n$.
Since $I$ is a $T$-ideal, we have
$$\eta(z(x_{n+1},\dots,x_{2n}) a(x_1,\dots,x_n))=z(x_{1},\dots,x_{n}) a((\mu(x_1),\dots,\mu(x_n))=z\mu(a)\in I.
$$
Since $z$ is an arbitrary element of $J$, the lemma is proved.
\endproof

Let $V$ and $U$ be ideals of $F\la X\ra$ and $V\subseteq U$. We introduce
the concept of the {\it separator} $V\div U$ as follows.

There exists
a unique smallest $T$-ideal $L$ of $F\la X\ra$ such that $V\subseteq UL$. (This is
an immediate corollary of the fir-property of $F\la X\ra$. Indeed, we can assume that $U\ne \{0\}$. Let
$(e_{\omega})_{\omega\in\Omega}$ be an $F\la X\ra$-basis of the free right ideal $U$,
and $L'$ be a $T$-ideal of $F\la X\ra$. Then
an element $\sum e_{\omega}f_{\omega}$ belongs to the product $UL'$ iff all the coefficients
$f_{\omega}$ belong to $L'$. Therefore $UL'\cap UL''\cap\dots=U(L'\cap L''\cap\dots)$
for arbitrary $T$-ideals $L', L'',\dots $) 
We denote this $T$-ideal $L$ by $V\div U$.  

\begin{rk} There also exists a least ideal $J$ such that $V\subseteq UJ$.
However $J$ is not necessarily a $T$-ideal even if both $U$ and $V$ are $T$-ideals of $F\la X\ra$.
Let for example, $U$ and $V$ be minimal $T$-ideals of $F\la X\ra$ containing $\sum_{\pi\in S_3} sign(\pi)x_{\pi(1)}x_{\pi(2)}x_{\pi(3)}$ and
$\sum_{\pi\in S_4} sign(\pi)x_{\pi(1)}x_{\pi(2)}x_{\pi(3)} x_{\pi(4)}$, respectively. Then
it is easy to see that the ideal $J$ consists of all polynomials with zero constant term. But $J$
is not a $T$-ideal since it is not invariant under the endomorphism $x_1\to 1,x_2\to 0,x_3\to 0,\dots$ of $F\la X\ra$.
\end{rk}

\begin{lemma}\label{var} Let $I$, $J$ and $U$ be $T$-ideals of $F\la X\ra$. Assume that $J\subseteq U$. Then
$$(UI:J)=(I:(J\div U)).$$ 
\end{lemma}

\proof We may assume that $U\ne \{0\}$. By Lemma \ref{incl}, for arbitrary ideals $V_1,V_2$ and $W$ of $F{\la X\ra}$, we have
  $UV_1W\subseteq UV_2$ iff $V_1W\subseteq V_2$. It follows that  $(I:(J\div U))=(UI:U(J\div U))\subseteq(UI:J)$ because $J\subseteq U(J\div U)$. 
  It remains to prove that $(UI:J)\subseteq(I:(J\div U))$.
 
   Let $a=a(x_1,\dots, x_n)\in (UI:J)$, and so $za\in UI$ for every $z\in J$. With respect to the basis $B=(e_{\alpha})$
of the (right) ideal $U$, given by Lemma \ref{basis}, we have $z=\sum e_{\alpha}f_{z,\alpha}$, where $f_{z,\alpha}\in F\la X\ra$.
Hence $f_{z,\alpha}a\in I$ for every $f_{z,\alpha}$.
  
  The ideal $L=J\div U$ is the minimal $T$-ideal of $F\la X\ra$ containing all the polynomials $f_{z,\alpha}=f_{z,\alpha}(x_1,\dots,x_n).$ 
  To complete the proof, it suffices to show that for every endomorphism $\mu: F\la X\ra \to F\la X\ra$ and every polynomial $g=g(x_1,\dots,x_n)$, we have  $(\mu(f_{z,\alpha})g)a\in I$.
  
  Recall that $\phi(x_i)=x_{i+1}$. Since $J$ is a $T$-ideal, we have $\phi^n(z)g=\sum \phi^n(e_{\alpha})\phi^n(f_{z,\alpha})g\in J$,
  and this is a basis decomposition of $\phi^n(z)g$ by Lemma \ref{basis}. Therefore $\phi^n(f_{z,\alpha})ga\in I$, that is,
  $f_{z,\alpha}(x_{n+1},\dots,x_{2n})g(x_1,\dots,x_n)a(x_1,\dots,x_n)\in I$. Since $I$ is a $T$-ideal, it is invariant
  under an endomorphism of $F\la X\ra$ such that $x_1\mapsto x_1,\dots, x_n\mapsto x_n, x_{n+1}\mapsto \mu(x_{1}),\dots, x_{2n}\mapsto\mu(x_n)$.
  Hence $f_{z,\alpha}(\mu(x_1),\dots,\mu(x_{n}))g(x_1,\dots,x_n)a(x_1,\dots,x_n)\in I$, as required. 
  \endproof

For the next lemma, we need the construction of triangular product of modules
invented by B.I.Plotkin  and U.E.Kal'julaid (\cite{PV},\cite{K}). Note that every $R$-module $P$
over an $F$-algebra $R$ is a module over $\bar R$, the image of $R$ in $End_F(P)$.

Let $R_1$ and $R_2$ be two
$F$-algebras, $M_1$ a (right) $R_1$-module, and $M_2$ an $R_2$-module.  Then the algebra 
$Q=\bar R_1\oplus\bar R_2$ canonically becomes a subalgebra of $End_F(M_1\oplus M_2)$.
Let $\Phi$ consists of the $F$-linear operators on $M_1\oplus M_2$ mapping $M_1$ to $\{0\}$ and $M_2$ to $M_1$.
Then $R=Q+\Phi$ is a subalgebra  of $End_F(M_1\oplus M_2)$, and so the vector space $M=M_1+M_2$ is a faithful 
$R$-module. It is called the {\it triangular product} $M_1 \nabla M_2$ of the modules $M_1$ and $M_2$. (Observe
that $M_1$ is an $R$-submodule, but $M_2$ is not.)

Denote by $T(M_1)$ the $T$-ideal of $ F\la X\ra$ consisting of all polynomials $f$ vanishing under every
homomorphism $F\la X\ra\to \bar R_1$. Similarly, one defines $T(M_2)$ and $T(M_1\nabla M_2)$. Then $T(M_1\nabla M_2)=T(M_2)T(M_1)$. (See \cite{K} and
Appendix in \cite{V}. 
Of course, the inclusion $\supseteq$ follows from the definitions.) 

\begin{lemma}\label{07} Let $I$, $J$, and $K$ be  $T$-ideals of $F\la X\ra$ and $K\not\subseteq I$.  
Then $$(IJ:K)=(I:K)J.$$  
\end{lemma} 
\proof The inclusion $\supseteq$ is obvious. To obtain $\subseteq$, we consider the right  $F\la X\ra$-modules
$M_1=F\la X\ra/J$,  $M_2=F\la X\ra/I$, and their triangular product $M=M_1\nabla M_2$ which is an $R$-module
as in the definition of triangular product. Since $I$ and $J$ are $T$-ideals, we have $T(M_1)=J$, $T(M_2)=I$, and therefore $T(M)=T(M_2)T(M_1) = IJ$. Also observe that $T(M_2K)=(I:K)$ by Lemma \ref{ann}.

For a $T$-ideal $L$ of $F\la X\ra$, we denote by $L_R$ the corresponding $T$-ideal of $R$, that is the minimal ideal of $R$ containing $\mu(L)$
for each of the homomorphisms $\mu: F\la X\ra\to R$. Note that $L_R$ is not only generated as an ideal, but in fact is spanned 
as an $F$-vector-space by all $\mu (f)$, where $f=f(x_1,\dots, x_n)\in L$. Indeed,
for arbitrary $r_1,r_2\in R$ the product $r_1\mu(f)r_2$ is the image of the polynomial $x_{n+1}f x_{n+2}$ under a homomorphism
$\eta$ such that $\eta(x_1)=\mu(x_1),\dots,\eta(x_n)=\mu(x_n),\eta(x_{n+1})=r_1,$ and $\eta(x_{n+2})=r_2$.

Since $K\not\subseteq I$, we have $M_2K\ne 0$, and so $N=MK_R\not\subseteq M_1$.
Now, using $F$-linear operators from $\Phi$, we observe that $N\supseteq M_1$.
Furthermore, $N$ regarded as $\bar R$-module, where $\bar R$ is the image of $R$ in $End(M_1+M_2K)$, is isomorphic to $M_1\nabla (M_2K)$.
Therefore $T(N)=T(M_2K)T(M_1)=(I:K)J$. To complete the proof, it suffices to show that $T(N)\supseteq (IJ:K)$.  

Assume that $a=a(x_1,\dots,x_n)\in (IJ:K)$, that is, $fa\in IJ$ for every $f=f(x_1,\dots,x_n)\in K$. Since $K$ is
a $T$-ideal, we also have $\bar f a \in IJ$, where $\bar f=f(x_{n+1},\dots, x_{2n}).$ Recall that $T(M)=IJ$, and so we have 
$v (\mu(\bar f)\mu(a))=0$ for every $v\in M$ and every homomorphism $\mu: F\la X\ra\to R$.
Since there are no variables $x_i$ involved in both $\bar f$ and $a$, we also have $(v (\mu_1(\bar f))\mu_2(a)=0$ for any
two endomorphisms $\mu_1,\mu_2: F\la X\ra\to R$,  and so $(M K_R) \mu_2(a)=0$ by the definition of $K_R$ because we chose
arbitrary $f\in K$. But $MK_R=N$, and therefore $a\in T(N)$, as desired.

\endproof

\section{Proofs of the theorems}
  
{\it Proof of Theorem \ref{theor}.} The remainder of the proof of Theorem \ref{theor} is similar to that in \cite{O} and \cite{OS}.
The 'if' part of the statement is obvious, and we prove the 'only if' part below.

If $V'=\{0\}$, we can choose $I_1=\dots=I_k=0$ and $I'_1=V'_1,\dots, I'_l=V'_l$. Hence we may suppose that $V'\ne \{0\}$
and induct on the number $l$ of the factors in the decomposition of $V'$. 

 If $V\supseteq V'_1$, then there is
an easy solution, namely, $I_1=V_1,\dots, I_{k-1}=V_{k-1},I_k=V_k V'_2\dots V'_l$ and
 $I'_1=V, I'_2=V'_2,\dots, I'_l=V'_l$.
So we may assume that  $V\not\supseteq V'_1$, and consequently, $l\ge 2$.

Thus, there is $j$ ($1\le j\le k$) such that $V_1\dots V_{j-1}\supseteq V'_1$ but $V_1\dots V_{j-1}V_j\not\supseteq V'_1$.
The $T$-ideal $U=V_1\dots V_{j-1}$ is not $\{0\}$ since $U\supseteq V'_1\supseteq V'\ne \{0\}$. Now the separator $ V'_1\div U=W_1$ is a $T$-ideal of $F\la X\ra$, and by Lemma \ref{var},
we get $(UV_j:V'_1)=(V_j:W_1)$. The ideal $W_2=(UV_j:V'_1)=(V_j:W_1)$ is a $T$-ideal by Lemma \ref{ann}. Recall that $W_1(V_j:W_1)\subseteq V_j$
by definition, whence

\begin{equation}\label{w1w2}
W_1W_2\subseteq V_j 
\end{equation}
Since $U=V_1\dots V_{j-1}$, it follows from the definition of $W_1$ that

 \begin{equation}\label{timesw1}
V_1\dots V_{j-1}W_1\supseteq V'_1
\end{equation}
Now we use Lemma \ref{07} with $I=V_1\dots V_j$, $J=V_{j+1}\dots V_k$, $K=V'_1$, whence

 \begin{equation} \label{from}
 (V: V'_1)=W_2V_{j+1}\dots V_k
\end{equation}
by the definition of $W_2$. It follows from the inclusion $V\supseteq V'=V'_1(V'_2\dots V'_l)$ 
that $V'_2\dots V'_l\subseteq (V: V'_1)$, and therefore by (\ref{from}), we have 

 \begin{equation}\label{minus1}
W_2V_{j+1}\dots V_k\supseteq V'_2\dots V'_l
\end{equation}

The inclusion (\ref{minus1}) and the inductive hypothesis imply the existence of
a $T$-ideal $L$ of $F\la X\ra$ with two $T$-ideal factorizations

 \begin{equation}\label{induc}
L=L_j L_{j+1}\dots L_k = L'_2\dots L'_l
\end{equation}
such that 

\begin{equation} \label{compar}
L_j\subseteq W_2, L_{j+1}\subseteq V_{j+1},\dots, L_k\subseteq V_k;\;\; V'_2\subseteq L'_2,\dots, V'_l\subseteq L'_l
\end{equation}

To complete the proof, we set 

 \begin{equation}\label{solution}
 I_1=V_1,\dots, I_{j-1}=V_{j-1}, I_j=W_1L_j, I_{j+1}=L_{j+1},\dots, I_k=L_k, 
 \end{equation}
and
 \begin{equation}\label{solution'}
I'_1=V_1\dots V_{j-1}W_1, I'_2=L'_2,\dots, I'_l=L'_l
\end{equation}

Then it follows from (\ref{induc}--\ref{solution'}) that $I_1\dots I_k =I'_1\dots I'_l$ $(=I)$.
Also $I_1\subseteq V_1,\dots, I_{j-1}\subseteq V_{j-1}, I_{j+1}\subseteq V_{j+1},\dots, I_k\subseteq V_k$ by (\ref{solution})
and (\ref{compar}), and $I_j=W_1L_j\subseteq W_1W_2\subseteq V_j$ by (\ref{compar}) and (\ref{w1w2}).
Finally, $V'_1\subseteq I'_1$ by (\ref{solution'}) and (\ref{timesw1}), and $V'_2\subseteq I'_2,\dots, V'_l\subseteq I'_l$
by (\ref{solution'})  and (\ref{compar}). The theorem is proved.

\medskip

\begin{rk}  As G.Bergman and J.Lewin proved for $F\la X\ra$ that both the monoids of non-zero ideals and of non-zero $T$-ideals are free,
it would be interesting to determine if one can omit all $T$-s in the formulations of Lemma \ref{07} and Theorem \ref{theor} and replace
the infinite set $X$ by a finite one with $card (X)\ge 2$. Also it would  be interesting to see to what extent the results can be generalized to more general firs.
\end{rk}

\medskip

{\it Proof of Theorem \ref{theor'}.} Let $G$ be the free group freely generated by $X=\{x_1,x_2,\dots\}$ and $FG$
the group algebra of $G$ over  a field $F$. Then the set of varieties of group representations over $F$ is in
Galois correspondence with the set of fully invariant ideals of $FG$ (see
\cite{PV}, Theorem I.2.1.2 or Section I.1 in \cite{V}). By definition, these ideals are invariant under the
endomorphisms of the free group $G$. Therefore it suffices to prove the analog of Theorem \ref{theor} for 
fully invariant ideals of $FG$. To obtain such a proof, we replace $F\la X\ra $ by $FG$,
$T$-ideals by  fully invariant ideals, and make the following minor alternation of our argument.

It is easy to see that $1$ and all the products $u=(1-x_{i_1}^{n_1})\dots(1-x_{i_t}^{n_t})$, where $t\ge 1, n_j\in \mathbb Z, n_j\ne 0,
i_j\ne i_{j+1},$ form  an $F$-basis $X^*$ of $FG$. Now we take 'monomials' $u$ of this form in the definition of Schreier set, replacing the word "prefix" by "left factor". By definition, $\deg u=|n_1|+\dots+|n_t|$. We suppose that $X^*$ is well ordered 
by degrees and lexicographically if the degrees of two monomials are equal ($1<1-x_1<1-x_1^{-1}<1-x_2<\dots$). 

If we now set $\psi(x_1)=1$ in the definition of the endomorphism $\psi$ (and $\psi(x_{i+1})=x_i$ for $i\ge 1$, as earlier), then
the proofs of lemmas \ref{schr} and \ref{basis} work for every fully invariant ideal $I$ of $FG$. The claim of Lemma \ref{incl} is also true since $FG$ is a fir by \cite{L}. 

There are no changes in the proof of the analogs of lemmas \ref{ann} and \ref{var}. In the definition of the triangular
product, now $R_1$ and $R_2$ are group algebras of groups $H_1$ and $H_2$, $\bar H_1$ and $\bar H_2$ are the canonical
images of $H_1$ and $H_2$ in the $F$-linear groups $GL(M_1)$ and $GL(M_2)$, and so $Q=\bar H_1\times\bar H_2$ is a subgroup of $GL(M_1\oplus M_2)$, $\Phi$ is the group of $F$-linear operators of $M_1\oplus M_2$
identical on the $F$-spaces $M_1$ and $(M_1\oplus M_2)/M_1$, and $R=FH$, where $H=Q\Phi$. In the definition of $T(M_1)$,
only those homomorphisms of group algebras are involved that are induced by group homomorpisms $G\to H_1$. The equality
$T(M_1\nabla M_2)=T(M_1)T(M_2)$ is Vovsi's theorem (\cite{V}, I.6.2). We must choose $r_1, r_2\in H$ in the proof of the analog of Lemma \ref{07}. 
Then the proof of Theorem \ref{theor} just turns into the proof of Theorem \ref{theor'}.

\medskip

{\bf Acknowledgments.} The author is very much obliged to G.Bergman who read the manuscript, shortened the proof of Lemma \ref{schr},
and suggested many small improvements.
The author is also grateful to Yu.A.Bahturin and A.L.Shmel'kin for useful discussion.

\begin{minipage}[t]{3 in}
\noindent Alexander Yu. Ol'shanskii\\ Department of Mathematics\\
Vanderbilt University \\ alexander.olshanskiy@vanderbilt.edu\\
 and\\ Department of
Higher Algebra\\ MEHMAT,
 Moscow State University\\

\end{minipage}

\end{document}